\documentclass[12pt, oneside, final, a4paper]{article}

\usepackage[T1]{fontenc}
\usepackage[utf8]{inputenc}
\usepackage[english]{babel}
\usepackage{indentfirst}

\usepackage{geometry}

\usepackage[tbtags]{amsmath}
\usepackage{amsfonts,amssymb,mathrsfs,amsthm}

\numberwithin{equation}{section}

\theoremstyle{plain}
\newtheorem{theorem}{Theorem}[section]
\newtheorem{lemma}[theorem]{Lemma}
\newtheorem{propos}[theorem]{Proposition}

\theoremstyle{definition}
\newtheorem{definition}[theorem]{Definition}
\newtheorem{remark}[theorem]{Remark}

\DeclareMathOperator{\Lip}{Lip}
\DeclareMathOperator{\dist}{dist}
\DeclareMathOperator{\supp}{supp}

\usepackage{csquotes} 
\usepackage[sorting=none, backend=biber]{biblatex}
\begin{document}

\title{Composition operators on Sobolev spaces with
a~metric measure space as a~domain.~II}
\author{Danil A. Sboev\thanks{%
    Department of Mechanics and Mathematics, Novosibirsk State
    University, 1 Pirogov St, Novosibirsk, 630090, Russian Federation. 
    E-mail: d.sboev@g.nsu.ru}
  \and
  Sergey K. Vodopyanov\thanks{%
    Sobolev Institute of Mathematics, Siberian Branch of
    the Russian Academy of Sciences, 4 Akademika Koptyuga Pr,
    Novosibirsk, 630090, Russian Federation.
    E-mail: vodopis@math.nsc.ru}}
\date{}

\maketitle

\bgroup
\let\thefootnote\relax
\egroup

\begin{abstract}
  We obtain a~description of the homeomorphisms
  which induce bounded composition operators
  on Sobolev spaces of functions
  on metric measure spaces.

  \textbf{Bibliography:} 8 items.

\end{abstract}

\tableofcontents

\section{Introduction}

This article continues the series begun in~\cite{VS25}.
The previous article investigated
necessary and sufficient conditions on a~homeomorphism
$\varphi \colon X \to Y$
between metric measure spaces
which yield the bounded composition operator
$$
\varphi^{*} \colon D^{1,p}(Y) \cap \Lip(Y) \to D^{1,p}(X), \quad 1
\le p < \infty,
$$
of Sobolev spaces
defined as
$u \mapsto \varphi^{*}(u) = u \circ \varphi$.

The main result of this article
is the following assertion.

\begin{theorem}\label{thm: main theorem}
  Consider two metric measure spaces
  $(X, d_{X}, \mu)$
  and
  $(Y, d_{Y}, \nu)$,
  where~$\mu$
  is a~Radon measure,
  as well as a~homeomorphism
  $\varphi \colon X \to Y$
  and two numbers
  $1 \le q < p < \infty$
  with
  $\frac{1}{\sigma} = \frac{1}{q} - \frac{1}{p}$.

  The mapping~$\varphi$
  induces the bounded composition operator
  $$
  \varphi^{*} \colon D^{1,p}(Y) \cap \Lip(Y) \to D^{1,q}(X)
  $$
  if and only if:

  {\rm(1)}
  $\varphi \in D^{1,q}_{\operatorname{loc}}(X; Y)$$;$

  {\rm(2)}
  $\varphi$
  has finite distortion$;$

  {\rm(3)}
  the operator distortion function
  $K_{q,p}(x, \varphi)$
  lies in
  $L^{\sigma}(X)$.

  Furthermore,
  we have the equality
  $\|K_{q,p}(\cdot, \varphi) \mid L^{\sigma}(X)\| = \|\varphi^{*}\|$.
\end{theorem}

\begin{remark}
  Theorem~\ref{thm: main theorem}
  in the case
  $p=q$
  was established~\cite{VS25} by the authors
  under weaker assumptions on~$X$ and~$Y$.
\end{remark}

In contrast to the previously available methods for proving necessity,
we obtain the regularity of the homeomorphism
and estimate the distortion function
resting on new ideas.

In this article,
the regularity follows from
our consideration of the additional composition operator
(the case
$1 \le q < p = \infty$):
\begin{equation}\label{eq: intro op comp lip}
  \varphi^{*} \colon \Lip(Y) \to D^{1,q}(X), \quad 1 \le q < \infty.
\end{equation}
The mappings of Sobolev class
$D^{1,q}(X;Y)$
clearly induce the composition operator \eqref{eq: intro op comp lip}.
This article establishes the converse.

\begin{theorem}
  Consider two metric measure spaces
  $(X, d_{X}, \mu)$
  and
  $(Y, d_{Y}, \nu)$.
  A~homeomorphism
  $\varphi \colon X \to Y$
  induces the bounded composition operator
  $$
  \varphi^{*} \colon \Lip(Y) \to D^{1,q}(X), \quad 1 \le q < \infty,
  $$
  if and only if
  $\varphi \in D^{1,q}(X; Y)$
  in the sense of Reshetnyak.
  Furthermore,
  $\|\varphi^{*}\|^{q} = \int\limits_{X}|D_{q}\varphi|^{q}\, d \mu$.
\end{theorem}

Note that
\cite{PV25} justifies a~similar result for the homeomorphisms
which map a~region of a~Carnot group into a~metric space.
In this article,
the proof rests additionally on some ideas and methods,
borrowed from the theory of
$K$-spaces.

The next step in studying necessary conditions on homeomorphisms
is to justify an~estimate for the operator distortion function
$K_{q,p}(x, \varphi)$.
Usually,
the integral of the distortion function
is bounded by an~expression involving some quasi-additive function~$\Phi$
defined on an~open subset of~$Y$.
Thanks to Lebesgue's theorem and
theorems about differentiating a~quasi-additive set function,
a~pointwise inequality for the distortion function
is established.

This article associates to the composition operator
$$
\varphi^{*} \colon D^{1,p}(Y) \cap \Lip(Y) \to D^{1,q}(X),
\quad 1 \le q < p < \infty,
$$
some finite regular Borel measure~$\widetilde{\Psi}$.
We show that
$\Phi$
is the trace of~$\widetilde{\Psi}$
on the open subsets of finite measure.
Therefore,
instead of the theorem about differentiating a~quasi-additive set function
we can apply the Radon--Nikodym theorem.

Since we do not require
$(X, d_{X}, \mu)$
and
$(Y, d_{Y}, \nu)$
to admit Vitali coverings,
Lebesgue's theorem may not yield pointwise estimates.
In this article we use the ``reverse'' H\"older's inequality.

\begin{propos}
  Consider three positive measurable functions~$f$,
  $g$,~and~$h$
  on some measure space
  $(Y, \Sigma, \nu)$
  and take two numbers
  $p$, $q \ge 1$
  with
  $\frac{1}{p} + \frac{1}{q} = 1$.

  The inequality
  $$
  \int\limits_{U} f(y) \, d \nu(y)
  \le
  \bigg(\int\limits_{U} g(y) \, d \nu(y) \bigg)^{\frac{1}{p}}
  \bigg(\int\limits_{U} h(y) \, d \nu(y) \bigg)^{\frac{1}{q}}
  $$
  holds for every
  $U \in \Sigma$
  if and only if
  $f(y) \le g(y)h(y)$
  holds $\nu$-a.e.
\end{propos}

We relegate the proof of this proposition to
Appendix~\ref{sec: reversed Holder}.

The structure of this article is as follows.

Section~2 is a~short collection of basic preliminary facts.

Section~3 deals with set functions
associated to composition operators.
We show that,
given a~mapping
which induces the bounded composition operator
$$
\varphi^{*} \colon D^{1,p}(Y) \cap \Lip(Y) \to D^{1,q}(X),
\quad 1 \le q < p \le \infty,
$$
we can construct some finite regular Borel measure~$\Psi$.

Section~4 contains a~description of
the Sobolev-class homeomorphisms in the sense of Reshetnyak.
Using the results of the theory of
$K$-spaces
(the criterion for~$o$-boundedness of a~family in a~$K$-space with
a~metric function)
and the results of Section~3,
we prove that
a~homeomorphism~$\varphi$
induces the bounded composition operator
$$
\varphi^{*} \colon \Lip(Y) \to D^{1,q}(X), \quad 1\le q < \infty,
$$
if and only if
$\varphi \in D^{1,q}(X;Y)$
in the sense of Reshetnyak.

Section~5 is devoted to proving the main result of this article.

Appendix~\ref{sec: reversed Holder}
contains a~proof of the ``reverse'' H\"older's inequality.

To make our exposition more complete,
Appendix~\ref{sec: K lem}
contains a~proof of a~particular case of the~$o$-boundedness criterion
for a~family of functions in
$L^{1}(X, \mu)$.


\section{Preliminaries}

The main properties of metric measure spaces,
Sobolev spaces,
and Sobolev mappings in the sense of Reshetnyak
are presented in \cite{VS25}.
Let us include here only some propositions
directly required in this article.

Recall that
a~metric measure space is a~triple
$(X, d_{X}, \mu)$,
where
$(X, d_{X})$
is a~separable metric space
and~$\mu$
is a~nontrivial regular Borel measure
which is finite on bounded subsets.

\subsection{Radon measures}

Say that
a~measure~$\mu$
is a~\emph{Radon measure}
on a~metric space
$(X, d_{X})$
whenever~$\mu$
is a~regular Borel measure,
$\mu(K) < \infty$
for every compact subset
$K \subset X$,
and every~$\mu$-measurable set~$A$
satisfies 
$$
\mu(A) = \sup\{ \mu(K) \mid A \supset K\text{ is a~compact subset} \}.
$$

Proposition 3.3.46 of \cite{HKNT15}
gives a~criterion for~$\mu$
to be a~Radon measure on a~metric measure space.
We can state a~particular case of that proposition as follows:
if
$(X, d_{X}, \mu)$
is a~metric measure space
and the space
$(X, d_{X})$
is locally compact
then~$\mu$
is a~Radon measure.

The main properties of Radon measures on Hausdorff locally compact spaces
are explained in \cite[Ch.~2]{Fed}.

\subsection{Pavlov's lemma}

The following lemma is originally proved as \cite[Lemma 4.2]{PV25}
for the measurable mappings defined on domains in Carnot groups
with values in a~metric space.
Its generalization to metric measure spaces
is presented in~\cite{VS25}.

Given an~open subset
$U \subset Y$
of a~metric space
$(Y, d_{Y})$,
denote by
$\overset{\circ}{\operatorname{Lip}}(U)$
the space of Lipschitz functions~$u$
with
$$\operatorname{dist}_{Y}(\supp u, Y \setminus U) > 0.$$

\begin{lemma}
  \label{lem: Pavlov}
  Consider a~metric measure space
  $(X, d_{X}, \mu)$,
  a metric space
  $(Y, d_{Y})$,
  and take
  $1 \le q < p = \infty$.
  Consider also a~measurable mapping
  $\varphi \colon X \to Y$,
  an~open subset
  $U \subset Y$,
  and some constant
  $A > 0$.

  Suppose that
  for every function
  $u \in \overset{\circ}{\operatorname{Lip}}(U)$
  the composition
  $u \circ \varphi$
  lies in the space
  $D^{1,q}(X)$
  and satisfies 
  \begin{equation}\label{eq: ineq Pavl}
    \bigg(
      \int\limits_{\varphi^{-1}(U)} | \nabla_{q} (u\circ \varphi)|^{q}\, d \mu
    \bigg)^{\frac{1}{q}}
    \le
    A \operatorname{Lip}(u ; U),
  \end{equation}
  where
  $\operatorname{Lip}(u; U) =
  \sup\big\{\frac{|u(x)-u(y)|}{d_Y(x,y)} \mid x\neq y, \ x,y\in U\big\}$.

  Then 
  \eqref{eq: ineq Pavl}
  holds for every function
  $u \in \operatorname{Lip}(Y)$
  as soon as
  $u \circ \varphi \in D^{1,q}(X)$.
\end{lemma}

\subsection{The area formula for homeomorphisms}

Given a~homeomorphism
$\varphi \colon X \to Y$,
there exist some locally integrable Borel function
$\mathcal{J} \varphi$
and some Borel measure
$\Phi^{s}$,
singular with respect to~$\mu$,
such that
$$
\nu(\varphi(B))
=
\int\limits_{B} \mathcal{J} \varphi d \mu + \Phi^{s}(B)
$$
for every Borel set
$B \subset X$.

Therefore,
we have the decomposition
$X = Z_{\varphi} \cup \Sigma_{\varphi} \cup X_{\rm r}$,
where
$Z_{\varphi}$
is the zero set of the Jacobian
$\mathcal{J}\varphi$
and
$\Sigma_{\varphi}$
is the singularity set of~$\varphi$
(the set on which the singular measure
  $\Phi^{s}$
is concentrated),
while
$X_r = X \setminus (Z_{\varphi}\cup \Sigma_{\varphi})$.
Observe that
$\varphi$
outside the set
$\Sigma_{\varphi}$
has Luzin's
$\mathcal{N}$-property
because
$\Phi^{s} = 0$
outside
$\Sigma_{\varphi}$.
For the inverse mapping
$\varphi^{-1}$
we can write the similar decomposition
$Y = Z_{\varphi^{-1}} \cup \Sigma_{\varphi^{-1}} \cup Y_{\rm r}$.
Removing Borel subsets of~$\nu$-measure zero from
$Z_{\varphi^{-1}}$
and
$\Sigma_{\varphi^{-1}}$,
we may assume that
both spaces are decomposed into Borel sets:
$X = Z_{\varphi} \cup \Sigma_{\varphi} \cup X_{\rm r}$
and
$Y = Z_{\varphi^{-1}} \cup \Sigma_{\varphi^{-1}} \cup Y_{\rm r}$;
furthermore,
$\varphi(\Sigma_{\varphi}) = Z_{\varphi^{-1}}$
and
$\varphi(Z_{\varphi}) = \Sigma_{\varphi^{-1}}$.
For more details,
see \cite[\S\S~2.8--2.10]{VS25}

\subsection{Sobolev mappings in the sense of Reshetnyak}

The following definition in the case of mappings
from a~Euclidean space into a~metric space
was originally introduced in~\cite{Resh97}.

\begin{definition}
  Given a~metric measure space
  $(X, d_X, \mu)$
  and a~metric space
  $(Y, d_{Y})$,
  a~measurable mapping
  $\varphi \colon X \to Y$
  lies in the \emph{Sobolev class
    $D^{1,p}(X; Y)$
  in the sense of Reshetnyak},
  where
  $1 \le p < \infty$,
  whenever there exists a~function
  $w \in L^{p}(X)$
  such that
  for every function
  $u_{z}(y) = d(y,z)$
  of
  $z \in Y$
  the composition
  $[\varphi]_{z} = u_{z} \circ \varphi$
  belongs to 
  $D^{1,p}(X)$
  and the relation
  \begin{equation}\label{eq: def of Sob map}
    |\nabla_{p}[\varphi]_{z}| \le w \quad \text{$\mu$-a.e.\ in~$X$}
  \end{equation}
  holds pointwise.

  A~measurable mapping~$\varphi$
  lies in the \emph{Sobolev class}
  $D^{1,p}_{\operatorname{loc}}(X;Y)$
  whenever each point
  $x \in X$
  has a~neighborhood
  $U_{x}$
  in~$X$
  such that
  $\varphi \in D^{1,p}(U_{x}; Y)$.
\end{definition}

Among all the majorants~$w$
satisfying 
\eqref{eq: def of Sob map},
the smallest one exists,
which we denote by
$|D_{p} \varphi|$.

\subsection{Estimate for 
the gradient of the composition}

We also need the following statement,
whose proof appeared in \cite{VS25}.

\begin{propos}
  \label{prop: chain rule}
  If
  $u \colon Y \to \mathbb{R}$
  is a~Lipschitz function
  and
  $\varphi \colon X \to Y$
  is a~measurable mapping of Sobolev class
  $D^{1,p}_{\operatorname{loc}}(X; Y)$
  then the composition
  $u \circ \varphi$
  is a~function of class
  $D^{1,p}_{\operatorname{loc}}(X)$
  and 
  $$
  |\nabla_{p}(u\circ \varphi)|(x)
  \le
  |\nabla_{p} u|(\varphi(x)) |D_{p}\varphi|(x)
  \quad \text{for $\mu$-a.e.\ }\quad x \in X.
  $$
\end{propos}

\subsection{Mappings with finite distortion.
The distortion function}

Consider a~homeomorphism
$\varphi \colon X \to Y$
of Sobolev class
$D^{1,q}_{\operatorname{loc}}(X; Y)$
with
$1 \le q < \infty$.

\begin{definition}
  The homeomorphism~$\varphi$
  is called a~\emph{mapping with finite distortion}
  whenever
  $$
  |D_{q}\varphi|(x) = 0 \quad \text{for $\mu$-a.e.\ }
  \quad x \in Z_{\varphi} = \{x' \in X \mid \mathcal{J} \varphi(x')=0\}.
  $$
\end{definition}

Introduce also the operator distortion function.

\begin{definition}
  Consider a~homeomorphism
  $\varphi \colon X \to Y$
  of class
  $D^{1,q}_{\operatorname{loc}}(X;Y)$
  with finite distortion
  and
  $1\le q \le p < \infty$.
  The \emph{operator distortion function}
  $K_{q,p}(\cdot, \varphi)$
  is defined as
  $$
  K_{q,p}(x, \varphi) =
  \begin{cases}
    \frac{|D_{q}\varphi|(x)}{(\mathcal{J}\varphi(x))^{\frac{1}{p}}}
    &\text{if } \mathcal{J}\varphi(x) \neq 0,\\
    0
    &\text{otherwise.}
  \end{cases}
  $$
\end{definition}

Since the distortion of~$\varphi$
is finite,
the function
$K_{q,p}(\cdot, \varphi)$
is well-defined.

The symbol ``\qedsymbol'' indicates the end of a~proof.
Some theorems consist of a~sequence of proofs of several claims,
with the same symbol indicating the end of the proof of each claim.


\section{Set functions and composition operators}\label{par: set func}

Consider a~metric measure space
$(X, d_{X}, \mu)$
as well as a~measurable mapping
$\varphi \colon X \to Y$.
Assume that~$\varphi$
induces the bounded composition operator
\begin{equation}\label{eq: comp op}
  \varphi^{*}
  \colon
  D^{1,p}(Y) \cap \Lip(Y)
  \to
  D^{1,q}(X) \quad (1 \le q < p \le \infty),
\end{equation}
and so
$u\circ \varphi \in D^{1,q}(X)$
when
$u \in D^{1,p}(Y) \cap \Lip(Y)$
and we have 
$$
\bigg(\int\limits_{X} |\nabla_q (u \circ \varphi)|^{q}\, d
\mu\bigg)^{\frac{1}{q}}
\le \|\varphi^{*}\|
\bigg(\int\limits_{Y} |\nabla_p u|^{p}\, d \nu \bigg)^{\frac{1}{p}}
$$
for
$1 \le q < p < \infty$,
while if
$1 \le q < p = \infty$
then
$$
\bigg(\int\limits_{X} |\nabla_q (u \circ \varphi)|^{q}\, d
\mu\bigg)^{\frac{1}{q}}
\le \|\varphi^{*}\|
\Lip(u),
$$
where
$\Lip(u) = \Lip(u; Y)$.

In the case
$1 \le q < p < \infty$
it is known that
a~quasi-additive set function is related to
the composition operator.
Let us introduce the necessary concepts.

Denote by~$\mathscr{O}$
some system of open sets closed under finite unions.

\begin{definition}
  Say that
  $\Phi \colon \mathscr{O} \to [0, \infty]$
  is a~\emph{quasi-additive set function}
  whenever
  for all
  $n \in \mathbb{N}$
  and
  $U_1$,
  \dots,
  $U_n$,
  $U \in \mathscr{O}$
  the inclusion
  $\bigcup\limits_{k=1}^{n} U_n \subset U$
  implies that
  \begin{equation}\label{eq: qa def}
    \sum\limits_{k=1}^{n} \Phi(U_k) \le \Phi(U).
  \end{equation}
\end{definition}

If this becomes an~equality
for every finite disjoint tuple
$\{U_k \in \mathscr{O} \mid k=1,\dots,n\}$
and
$U = \bigcup\limits_{k=1}^{n} U_k$
then the function~$\Phi$
is called \emph{finitely additive}.
If~$\mathscr{O}$
is closed under countable unions
and the equality in~\eqref{eq: qa def}
holds for every countable collection
then~$\Phi$
is called \emph{countably additive}.

Suppose that~$\varphi$
induces the composition operator~\eqref{eq: comp op}
for
$1 \le q < p < \infty$.
Fix an~open subset
$U \subset Y$.
Consider the restriction
$$
\varphi^{*}\vert U \colon D^{1,p}(U) \cap \overset{\circ}{\Lip}(U)
\to D^{1,q}(X)
$$
of the composition operator.
Given an~open subset
$U \subset Y$,
put
\begin{equation}\label{eq: qa comp def}
  \Phi(U) = \|\varphi^{*}\vert U\|^{\sigma}, \quad
  \text{where}\quad \frac{1}{\sigma} = \frac{1}{q} - \frac{1}{p}.
\end{equation}

\begin{propos}\label{prp: qa comp op}
  Consider two metric measure spaces
  $(X, d_{X}, \mu)$
  and
  $(Y, d_{Y}, \nu)$
  and a~measurable mapping
  $\varphi \colon X \to Y$
  which induces the bounded composition operator~\eqref{eq: comp op}
  for
  $1 \le q < p < \infty$.

  The set function
  $\Phi(U)$
  defined in 
  \eqref{eq: qa comp def}
  on the open subsets
  $U \subset X$
  is monotone and countably additive.
\end{propos}

Proposition~\ref{prp: qa comp op}
is proved in \cite[Lemma 3.1]{VE22}
in the case of Carnot groups;
it does not rely on~$\varphi$
being a~homeomorphism.
Cf.\ the first proof in \cite{VU98}.
The argument of~\cite{VE22}
carries over to the case of metric measure spaces with obvious modifications.

Consider another set function
associated to the composition operator.
Given a~Borel subset
$B \subset Y$,
define the set function~$\Psi$
as
\begin{equation}\label{eq: meas def}
  \Psi(B)^{\frac{1}{\sigma}}
  = \sup \bigg(
    \int\limits_{\varphi^{-1}(B)} |\nabla_{q} (u\circ \varphi)|^{q}\, d \mu
  \bigg)^{\frac{1}{q}},
\end{equation}
where the supremum is over all functions
$u \in \Lip(Y)$
with 
$$
\int\limits_{B} |\nabla_p u |^{p}\, d \nu \le 1.
$$

The first question that arises
is whether the functions~$\Phi$ and~$\Psi$
coincide on open sets.
It is obvious that
$\Phi(U) \le \Psi(U)$
for every open subset
$U \subset Y$.
Let us verify the reverse inequality.

\begin{propos}\label{prp: funct eq}
  Consider two metric measure spaces
  $(X, d_{X}, \mu)$
  and
  $(Y, d_{Y}, \nu)$,
  as well as a~measurable mapping
  $\varphi \colon X \to Y$
  which induces the bounded composition operator \eqref{eq: comp op}
  for
  $1 \le q < p < \infty$.

  Then
  $\Phi(U) = \Psi(U)$
  for every open subset
  $U \subset Y$
  of finite measure,
  $\nu(U) < \infty$.
\end{propos}

\begin{proof}
  Fix some open subset~$U$
  of~$Y$
  and put
  $$
  [U]_{- \delta} = \{y \in Y \mid \dist_{Y}(y, Y \setminus U) > \delta\}.
  $$
  Then
  $U = \bigcup\limits_{j \in \mathbb{N}} [U]_{- \frac{1}{j}}$.

  Fix arbitrary
  $\varepsilon$,
  $\varepsilon' > 0$.
  It is clear that
  in the definition of the functions~$\Phi$ and~$\Psi$
  we can take the supremum over the bounded functions~$u$.
  Suppose that
  $u \in \Lip(Y) \cap D^{1,p}(U)$
  with
  $|u| \le M$.
  For sufficiently large~$j$
  we have 
  $$
  \int\limits_{\varphi^{-1}(U)}
  |\nabla_{q} (u\circ \varphi)|^{q}\, d \mu
  - \varepsilon
  <
  \int\limits_{\varphi^{-1}([U]_{-\frac{1}{j}})}
  |\nabla_{q} (u\circ \varphi)|^{q}\, d \mu.
  $$
  On the other hand,
  for sufficiently large~$j$
  we have 
  $l^{p} \nu(U \setminus [U]_{-\frac{1}{j}}) < \varepsilon'$,
  where
  $l = \Lip(u)$.

  By analogy with the proof of Lemma~\ref{lem: Pavlov},
  for
  $K > 0$
  introduce the function
  $R_{K} \colon \mathbb{R} \to \mathbb{R}$
  as
  $$
  R_{K}(x) =
  \begin{cases}
    K-x
    &\text{if }
    x > \frac{K}{2}, \\
    x
    &\text{if }
    -\frac{K}{2} \le x \le \frac{K}{2},\\
    -K-x
    &\text{if }
    x < - \frac{K}{2}.
  \end{cases}
  $$
  Choose
  $k \in \mathbb{N}$
  such that
  $$
  |R_{M / 2^{k}} \circ R_{M / 2^{k-1}} \circ \cdots \circ R_{M} \circ u|
  \le \frac{l}{2j}.
  $$
  Denote the resulting function by~$v$.
  Extend~$v$
  to a~function~$\overline{u}$
  satisfying
  $\overline{u}(y) = v(y)$
  for
  $y \in Y \setminus [U]_{-\frac{1}{j}}$
  and
  $\overline{u}(y) = 0$
  for
  $y \in Y \setminus [U]_{- \frac{1}{2j}}$.
  Then for sufficiently large~$j$
  we have 
  \begin{multline*}
    \int\limits_{\varphi^{-1}(U)}
    |\nabla_{q} (u\circ \varphi)|^{q}\, d \mu
    - \varepsilon
    <
    \int\limits_{\varphi^{-1}([U]_{-\frac{1}{j}})}
    |\nabla_{q} (u\circ \varphi)|^{q}\, d \mu
    \\
    =
    \int\limits_{\varphi^{-1}([U]_{-\frac{1}{j}})}
    |\nabla_{q} (\overline{u} \circ \varphi)|^{q}\, d \mu
    \le
    \int\limits_{\varphi^{-1}(U)}
    |\nabla_{q} (\overline{u} \circ \varphi)|^{q}\, d \mu
    \\
    \le
    \Phi(U)^{\frac{q}{\sigma}}
    \bigg(
      \int\limits_{U} |\nabla_{p} \overline{u}|^{p}\, d \nu
    \bigg)^{\frac{q}{p}}
    \le
    \Phi(U)^{\frac{q}{\sigma}}
    \bigg(
      \int\limits_{[U]_{-\frac{1}{j}}} |\nabla_{p} \overline{u}|^{p}\, d \nu
      + \varepsilon
    \bigg)^{\frac{q}{p}} \\
    \le
    \Phi(U)^{\frac{q}{\sigma}}
    \bigg(
      \int\limits_{U} |\nabla_{p} u|^{p}\, d \nu
      + \varepsilon'
    \bigg)^{\frac{q}{p}}.
  \end{multline*}
  Since
  $\varepsilon$,
  $\varepsilon' > 0$
  are arbitrary,
  it follow that
  $$
  \int\limits_{\varphi^{-1}(U)}
  |\nabla_{q} (u\circ \varphi)|^{q}\, d \mu
  \le
  \Phi(U)^{\frac{q}{\sigma}}
  \bigg(
    \int\limits_{U} |\nabla_{p} u|^{p}\, d \nu
  \bigg)^{\frac{q}{p}}
  $$
  for an~arbitrary bounded function
  $u \in \Lip(Y) \cap D^{1,p}(U)$.
  From this we infer that
  $\Psi(U) \le \Phi(U)$.
\end{proof}

The second natural question is
whether the function $\Psi$
is a~measure.

\begin{theorem}\label{thm: comp op measure}
  Consider two metric measure spaces
  $(X, d_{X}, \mu)$
  and
  $(Y, d_{Y}, \nu)$,
  where~$\mu$
  is a~Radon measure,
  as well as a~continuous mapping
  $\varphi \colon X \to Y$
  which induces the bounded composition operator~\eqref{eq: comp op}
  for
  $1 \le q < p < \infty$.

  Then the set function~$\Psi$
  defined as~\eqref{eq: meas def}
  extends to a~finite regular Borel measure
  $\widetilde{\Psi}$
  as
  $$
  \widetilde{\Psi}(A)
  =
  \inf\{ \Psi(U) \mid A \subset U\text{ is an open subset} \},
  \quad A \subset Y.
  $$
\end{theorem}

\begin{proof}
  We split the proof of Theorem~\ref{thm: comp op measure}
  into several propositions.
  \medskip

  \textbf{Proposition 1.}
  The set function~$\Psi$
  is subadditive on disjoint Borel sets.

  Take a~countable collection
  $(B_{j})_{j\in \mathbb{N}}$
  of disjoint Borel subsets
  and an~arbitrary function
  $u \in \Lip(Y)$
  with
  $|\nabla_p u | \in L^{p}(\bigcup\limits_{j\in \mathbb{N}} B_j)$.
  Applying H\"older's inequality
  (recall that
  $\frac{q}{\sigma} + \frac{q}{p} = 1$),
  we deduce the following relations:
  \begin{multline*}
    \int\limits_{\varphi^{-1}(\bigcup\limits_{j\in \mathbb{N}} B_{j})}
    |\nabla_q (u \circ \varphi)|^{q} \, d \mu
    =
    \sum\limits_{j\in \mathbb{N}}
    \int\limits_{\varphi^{-1}(B_j)}
    |\nabla_{q} (u \circ \varphi)|^{q} \, d \mu
    \\
    \le
    \sum\limits_{j \in \mathbb{N}}
    \Psi(B_{j})^{\frac{q}{\sigma}}
    \bigg(
      \int\limits_{B_j} |\nabla_{p} u|^{p}\, d \nu
    \bigg)^{\frac{q}{p}}
    \le
    \bigg(
      \sum\limits_{j\in \mathbb{N}} \Psi(B_j)
    \bigg)^{\frac{q}{\sigma}}
    \bigg(
      \sum\limits_{j \in \mathbb{N}}
      \int\limits_{B_j}
      |\nabla_{p} u|^{p}\, d \nu
    \bigg)^{\frac{q}{p}}\\
    =
    \bigg(
      \sum\limits_{j\in \mathbb{N}} \Psi(B_j)
    \bigg)^{\frac{q}{\sigma}}
    \bigg(
      \int\limits_{\bigcup\limits_{j\in \mathbb{N}} B_j}
      |\nabla_{p} u|^{p}\, d \nu
    \bigg)^{\frac{q}{p}}.
  \end{multline*}
  From this we infer that
  $\Psi(\bigcup\limits_{j\in \mathbb{N}} B_j)
  \le \sum\limits_{j\in \mathbb{N}} \Psi(B_j)$.
  \qed
  \medskip

  \textbf{Proposition 2.}
  If~$U$
  is a~Borel subset
  and~$W$
  is an~open subset with
  $U \subset W$
  then
  $\Psi(U) \le \Psi(W)$.

  Let us prove the following claim:
  if~$K$
  is a~compact subset of~$Y$
  and~$W$
  is an~open subset with
  $K \subset W$
  then
  $\Psi(K) \le \Psi(W)$.

  Take an~arbitrary function
  $u \in \Lip(Y) \cap D^{1,p}(K)$
  and fix some number
  $\varepsilon > 0$.

  Since~$K$
  is compact,
  we may assume that
  the function~$u$
  is bounded:
  $|u| \le M$.
  Put
  $\Lip(u) = l$
  and
  $[U]_{\gamma} = \{y \in Y \mid \dist_{Y}(y, U) < \gamma \}$,
  where
  $\gamma > 0$.
  For
  $K > 0$
  denote by
  $R_{K}$
  the function as in the proof of Proposition~\ref{prp: funct eq}.

  Since~$K$
  is a~compact set,
  while~$W$
  is an~open set,
  there exists
  $\delta > 0$
  such that
  $\dist_{Y}(K, Y \setminus W) = \delta$.
  Choose now
  $\gamma < \delta$
  such that
  $$
  l^{p} \cdot \nu([K]_{\gamma} \setminus K) \le \varepsilon.
  $$
  Choose also
  $k \in \mathbb{N}$
  such that
  $$
  |R_{M / 2^{k}} \circ R_{M / 2^{k-1}} \circ \cdots \circ R_{M} \circ u|
  \le \gamma.
  $$
  Denote the resulting function by~$v$.
  Extend~$v$
  to a~function
  $\overline{u}$
  satisfying
  $\overline{u}(y) = v(y)$
  for
  $y \in K$
  and
  $\overline{u}(y) = 0$
  for
  $y \notin Y \setminus [K]_{\gamma}$.
  The construction implies that 
  \begin{align*}
    \int\limits_{K} |\nabla_{p} u|^{p}\, d \nu
    & =
    \int\limits_{K} |\nabla_{p} \overline{u}|^{p}\, d \nu, \\
    \int\limits_{W} |\nabla_{p} \overline{u}|^{p}\, d \nu
    & \le
    \int\limits_{K} |\nabla_{p} \overline{u}|^{p}\, d \nu + \varepsilon.
  \end{align*}
  We also have the equality
  \begin{align*}
    \int\limits_{\varphi^{-1}(K)}
    |\nabla_{q} (u\circ \varphi)|^{q}\, d \mu
    & =
    \int\limits_{\varphi^{-1}(K)}
    |\nabla_{q} (\overline{u}\circ \varphi)|^{q}\, d \mu.
  \end{align*}

  From that we deduce the following inequalities:
  \begin{multline*}
    \int\limits_{\varphi^{-1}(K)}
    |\nabla_{q} (u\circ \varphi)|^{q}\, d \mu
    =
    \int\limits_{\varphi^{-1}(K)}
    |\nabla_{q} (\overline{u} \circ \varphi)|^{q} \, d \mu
    \\
    \le
    \int\limits_{\varphi^{-1}(W)}
    |\nabla_{q} (\overline{u} \circ \varphi)|^{q} \, d \mu
    \le
    \Psi(W)^{\frac{q}{\sigma}}
    \bigg(
      \int\limits_{W} |\nabla_{p} \overline{u}|^{p}\, d \nu
    \bigg)^{\frac{q}{p}}
    \\
    \le
    \Psi(W)^{\frac{q}{\sigma}}
    \bigg(
      \int\limits_{K} |\nabla_{p} \overline{u}|^{p}\, d \nu + \varepsilon
    \bigg)^{\frac{q}{p}}
    =
    \Psi(W)^{\frac{q}{\sigma}}
    \bigg(
      \int\limits_{K} |\nabla_{p} u|^{p}\, d \nu + \varepsilon
    \bigg)^{\frac{q}{p}}.
  \end{multline*}
  Since
  $\varepsilon > 0$
  and the function
  $u \in \Lip(Y) \cap D^{1,p}(K)$
  are arbitrary,
  we conclude that
  $\Psi(K) \le \Psi(W)$.

  Consider now
  an~arbitrary Borel subset~$B$
  of
  $Y$,
  an~open subset~$W$
  with
  $B \subset W$,
  as well as an~arbitrary compact subset~$K$
  of
  $\varphi^{-1}(B)$.
  Then
  $\varphi(K)$
  is a~compact subset of~$B$.
  We find that 
  $$
  \int\limits_{K} |\nabla_{q} (u\circ \varphi)|^{q}\, d \mu
  \le
  \Psi(\varphi(K))^{\frac{q}{\sigma}}
  \bigg(
    \int\limits_{\varphi(K)}
    |\nabla_p u|^{p}\, d \nu
  \bigg)^{\frac{q}{p}}
  \le
  \Psi(W)^{\frac{q}{\sigma}}
  \bigg(
    \int\limits_{B}
    |\nabla_p u|^{p}\, d \nu
  \bigg)^{\frac{q}{p}}.
  $$
  Since~$\mu$
  is a~Radon measure
  and the compact subset
  $K \subset \varphi^{-1}(B)$
  is arbitrary,
  we infer that
  $$
  \int\limits_{\varphi^{-1}(B)} |\nabla_{q} (u\circ \varphi)|^{q}\, d \mu
  \le
  \Psi(W)^{\frac{q}{\sigma}}
  \bigg(
    \int\limits_{B}
    |\nabla_p u|^{p}\, d \nu
  \bigg)^{\frac{q}{p}}.
  $$
  Thus,
  $\Psi(B) \le \Psi(W)$.
  \qed
  \medskip

  \textbf{Proposition 3.}
  The function~$\Psi$
  is subadditive on open subsets.

  Take a~countable collection
  $(W_j)_{j \in \mathbb{N}}$
  of open subsets
  and put
  $C_1 = W_1$,
  then
  $C_2 = W_2 \setminus W_1$,
  $C_3 = W_3 \setminus (C_1 \cup C_2)$,
  and so forth.
  Then
  $C_j$
  is a~Borel subset of~$W_j$
  for all
  $j\in \mathbb{N}$
  and
  $\bigcup\limits_{j\in \mathbb{N}} C_j = \bigcup\limits_{j\in \mathbb{N}} W_j$;
  furthermore,
  $C_j \cap C_l = \varnothing$
  for
  $j \neq l$.
  Propositions~1 and~2 yield
  $$
  \Psi(\bigcup\limits_{j \in \mathbb{N}} W_j)
  =
  \Psi(\bigcup\limits_{j \in \mathbb{N}} C_j)
  \le
  \sum\limits_{j\in \mathbb{N}} \Psi(C_j)
  \le
  \sum\limits_{j \in \mathbb{N}} \Psi(W_j).
  $$
  Hence,
  $\Psi(\bigcup\limits_{j\in \mathbb{N}} W_j) \le \sum\limits_{j\in
  \mathbb{N}} \Psi(W_j)$
  for the open sets
  $(W_{j})_{j \in \mathbb{N}}$.
  \qed

  Define the function~$\widetilde{\Psi}$
  as
  $$
  \widetilde{\Psi}(A) = \inf \{ \Psi(W) \mid A \subset W
  \text{ is an open set} \},
  $$
  where~$A$
  is an~arbitrary subset of~$Y$.
  \medskip

  \textbf{Proposition 4.}
  $\widetilde{\Psi}$
  is a~finite regular Borel measure on~$Y$.

  The function
  $\widetilde{\Psi}$
  has the following obvious properties:

  {\rm(1)}
  $\widetilde{\Psi}$
  is a~monotone function;

  {\rm(2)}
  $\widetilde{\Psi}$
  is an~outer measure on~$Y$;

  {\rm(3)}
  $\widetilde{\Psi}$
  coincides with~$\Psi$
  on open subsets,
  in~particular,
  $\widetilde{\Psi}(Y) = \Psi(Y) = \|\varphi^{*}\|$;

  {\rm(4)}
  for every
  $A \subset Y$
  there exists a~$G_{\delta}$-set%
  \footnote{%
    Say that~$B$
    is a~$G_{\delta}$-set
    whenever
    $B = \bigcap\limits_{n\in \mathbb{N}} O_n$,
    where
    $O_n$
    is an~open subset for all
    $n \in \mathbb{N}$.
  }
  $B$
  with
  $\widetilde{\Psi}(A) = \widetilde{\Psi}(B)$.

  Let us show that
  $\widetilde{\Psi}$
  is a~Borel measure.
  To this end,
  use Carath\'eodory's criterion:
  take two subsets~$A$ and~$B$
  of~$Y$
  with
  $\dist_{Y}(A, B) = 2\delta > 0$.
  It suffices to show that
  $\widetilde{\Psi}(A \cup B) \ge \widetilde{\Psi}(A) + \widetilde{\Psi}(B)$.

  Fix an~arbitrary
  $\varepsilon > 0$.
  There exist open sets~$\widetilde{A}$ and~$\widetilde{B}$
  with
  $A \subset \widetilde{A}$
  and
  $B \subset \widetilde{B}$
  such that
  $\dist_{Y} (\widetilde{A}, \widetilde{B}) = \delta$
  and 
  $\widetilde{\Psi}(\widetilde{A} \cup \widetilde{B})
  \le \widetilde{\Psi} ( A \cup B) - \varepsilon$.

  Put
  $\widetilde{A}_n = \widetilde{A} \cap B_{Y}(o, n)$
  and
  $\widetilde{B}_n = \widetilde{B} \cap B_{Y}(o, n)$,
  where
  $o \in Y$
  is an~arbitrary point
  and
  $n \in \mathbb{N}$.
  It is obvious that
  $\nu(\widetilde{A}_n) < \infty$
  for every
  $n \in \mathbb{N}$,
  we have 
  $\widetilde{A} = \bigcup\limits_{n\in \mathbb{N}} \widetilde{A}_{n}$,
  and consequently,
  $\widetilde{\Psi}(\widetilde{A})
  = \lim\limits_{n\to \infty} \widetilde{\Psi}(\widetilde{A}_n)$;
  similar relations hold for~$\widetilde{B}$
  and
  $\widetilde{B}_n$.
  Since~$\Phi$ and~$\Psi$
  coincide on all open subsets of finite measure,
  see Proposition~\ref{prp: funct eq},
  and~$\Phi$
  is countably additive,
  we infer that
  \begin{multline*}
    \widetilde{\Psi}(A) + \widetilde{\Psi}(B)
    \le
    \widetilde{\Psi}(\widetilde{A}) + \widetilde{\Psi}(\widetilde{B}) \\
    =
    \lim\limits_{n \to \infty}
    (\widetilde{\Psi}(\widetilde{A}_n) + \widetilde{\Psi}(\widetilde{B}_n))
    =
    \lim\limits_{n \to \infty}
    (\Phi(\widetilde{A}_{n}) + \Phi(\widetilde{B}_n))
    =
    \lim\limits_{n \to \infty}
    \Phi(\widetilde{A}_{n} \cup \widetilde{B}_n) \\
    =
    \lim\limits_{n \to \infty}
    \widetilde{\Psi}(\widetilde{A}_{n} \cup \widetilde{B}_{n})
    = \widetilde{\Psi}(\widetilde{A} \cup \widetilde{B})
    \le \widetilde{\Psi}(A \cup B) - \varepsilon.
  \end{multline*}
  Since
  $\varepsilon > 0$
  is arbitrary,
  the claim follows.
  \qed

  The proof of Theorem~\ref{thm: comp op measure} is complete.
\end{proof}

Consider now the set function associated to the composition operator
for
$p = \infty$.
Suppose that
a~measurable mapping
$\varphi \colon X \to Y$
induces the bounded composition operator~\eqref{eq: comp op}
for
$p = \infty$.
It turns out that
some regular Borel measure on~$Y$
is related to this composition operator.

\begin{theorem}\label{thm: measure with lip}
  Consider two metric measure spaces
  $(X, d_{X}, \mu)$
  and
  $(Y, d_{Y}, \nu)$,
  as well as a~measurable mapping
  $\varphi \colon X \to Y$
  which induces the bounded composition operator~\eqref{eq: comp op}
  for
  $1 \le q < p = \infty$.

  On the Borel~$\sigma$-algebra
  define the following set function:
  $$
  \mathscr{B}(Y) \ni U \mapsto
  \Psi(U) =
  \sup
  \bigg(
    \int\limits_{\varphi^{-1}(U)}
    |\nabla_{q}(u\circ \varphi)|^{q}
    \, d \mu
  \bigg)^{\frac{\sigma}{q}},
  $$
  where the supremum is over all functions
  $u \in \Lip(Y)$
  with
  $\Lip(u) \le 1$.

  Then the function~$\Psi$
  extends to the regular Borel measure~$\widetilde{\Psi}$
  on~$Y$
  as
  $$
  \widetilde{\Psi}(A) =
  \inf\{ \Psi(B) \mid A \subset B, \ B \in \mathscr{B}(Y) \},
  \quad A \subset Y.
  $$
\end{theorem}

\begin{proof}
  We split the proof of Theorem~\ref{thm: measure with lip}
  into several propositions.
  \medskip

  \textbf{Proposition 1.}
  The function~$\Psi$
  is countably subadditive on the Borel sets.

  Indeed,
  assume that
  $(B_j)_{j\in \mathbb{N}} \subset \mathscr{B}(Y)$
  and take
  $u \in \operatorname{1-Lip}(Y)$.
  This yields 
  $$
  \int\limits_{\varphi^{-1}\big(\bigcup\limits_{j=1}^{\infty} B_j \big)}
  |\nabla_q(u\circ \varphi)|^{q}\, d \mu
  \le
  \sum\limits_{j=1}^{\infty}
  \int\limits_{\varphi^{-1}(B_j)}
  |\nabla_q(u\circ \varphi)|^{q}\, d \mu
  \le
  \sum\limits_{j=1}^{\infty} \Psi(B_j).
  $$
  From this we infer that
  $\Psi\big( \bigcup\limits_{j=1}^{\infty} B_j \big) \le
  \sum\limits_{j=1}^\infty \Psi(B_j)$.
  \qed

  Extend~$\Psi$
  regularly to an~arbitrary subset of~$Y$
  as
  $$
  \widetilde{\Psi}(A) =
  \inf\{ \Psi(B) \mid A \subset B, \ B \in \mathscr{B}(Y) \}.
  $$

  \textbf{Proposition 2.}
  $\widetilde{\Psi}$
  is an~outer measure on~$Y$.

  Proposition~2 can be verified directly.
  \qed
  \medskip

  \textbf{Proposition 3.}
  $\widetilde{\Psi}$
  is a~finite regular Borel measure on~$Y$.

  Since the composition operator is bounded,
  it follows that~$\widetilde{\Psi}$
  is finite.

  To verify that
  $\widetilde{\Psi}$
  is a~Borel measure,
  we use Carath\'eodory's criterion.

  Take two subsets~$A$ an~$B$
  of~$Y$
  with
  $\dist_{Y}(A,B) = 2\delta > 0$.
  For arbitrary
  $\varepsilon > 0$
  there is a~Borel subset
  $W \supset A \cup B$
  with
  $\widetilde{\Psi}(A \cup B) - \varepsilon > \Psi(W)$.
  Consider the sets
  $\widetilde{A} = W \cap (A)_{\delta/4}$
  and
  $\widetilde{B} = W \cap (B)_{\delta/4}$,
  where
  $(A)_{\delta}$
  is the set
  $\{y \in Y \mid \operatorname{dist}_{Y}(y, A) < \delta\}$,
  meaning the~$\delta$-neighborhood of~$A$.
  It is clear that~$\widetilde{A}$ and~$\widetilde{B}$
  are Borel sets.

  Fix
  $\varepsilon' > 0$.
  There exist~$1$-Lipschitz functions~$u$ and~$v$
  with
  $|u| \le N$
  and
  $|v| \le N$
  for some
  $N > 0$
  satisfying 
  \begin{align*}
    \int\limits_{\varphi^{-1}(\widetilde{A})} |\nabla_{q} (u\circ
    \varphi)|^{q}\, d \mu
    & \ge \Psi(A) - \varepsilon', \\
    \int\limits_{\varphi^{-1}(\widetilde{B})} |\nabla_{q} (v\circ
    \varphi)|^{q}\, d \mu
    & \ge \Psi(B) - \varepsilon'.
  \end{align*}
  Using~$u$ and~$v$,
  construct some functions
  $\overline{u}$,
  $\overline{v} \in \operatorname{1-Lip}(Y)$
  with the following properties:

  {\rm(i)}
  $|\overline{u}| \le \frac{\delta}{2}$
  and
  $|\overline{v}| \le \frac{\delta}{2}$;

  {\rm(ii)}
  $\overline{u}$
  and
  $\overline{v}$
  vanish identically outside the sets
  $(\widetilde{A})_{\delta/4}$
  and
  $(\widetilde{B})_{\delta/4}$
  respectively;

  {\rm(iii)}
  for $\mu$-a.e.\ point
  $x \in \varphi^{-1}(\widetilde{A})$
  we have 
  $|\nabla_{q} (u \circ \varphi)| = |\nabla_{q}(\overline{u}\circ \varphi)|$
  and the similar equality for $\overline{v}$.

  The existence of~$\overline{u}$ and~$\overline{v}$
  follows from the construction used in Proposition~\ref{prp: funct eq}.

  Then
  $\overline{u}+\overline{v} \in \operatorname{1-Lip}(Y)$
  and
  \begin{multline*}
    \Psi(\widetilde{A} \cup \widetilde{B})
    \ge
    \int\limits_{\varphi^{-1}(\widetilde{A}\cup \widetilde{B})}
    |\nabla_{q} (\overline{u} + \overline{v}) \circ \varphi|^{q}
    \, d \mu
    \\
    =
    \int\limits_{\varphi^{-1}(\widetilde{A})} |\nabla_{q}
    (\overline{u} \circ \varphi)|^{q}\, d \mu
    +
    \int\limits_{\varphi^{-1}(\widetilde{B})} |\nabla_{q}
    (\overline{v} \circ \varphi)|^{q}\, d \mu
    \\
    =
    \int\limits_{\varphi^{-1}(\widetilde{A})} |\nabla_{q} (u \circ
    \varphi)|^{q}\, d \mu
    +
    \int\limits_{\varphi^{-1}(\widetilde{B})} |\nabla_{q} (v \circ
    \varphi)|^{q}\, d \mu
    \ge
    \Psi(\widetilde{A}) + \Psi(\widetilde{B}) - 2 \varepsilon'.
  \end{multline*}

  From this we deduce that
  \begin{multline*}
    \widetilde{\Psi}(A \cup B) - \varepsilon \ge \widetilde{\Psi}(W) = \Psi(W)
    \ge \Psi(\widetilde{A} \cup \widetilde{B}) \\
    \ge \Psi(\widetilde{A}) + \Psi(\widetilde{B}) - 2 \varepsilon'
    \ge \widetilde{\Psi}(A) + \widetilde{\Psi}(B) - 2 \varepsilon'.
  \end{multline*}
  Since
  $\varepsilon$,
  $\varepsilon' > 0$
  are arbitrary,
  we obtain the claim.
  The measure~$\widetilde{\Psi}$
  is regular by construction.
  \qed

  The proof of Theorem~\ref{thm: measure with lip} is complete.
\end{proof}

\begin{remark}
  In the case
  $1 \le q < p = \infty$
  we can also introduce the quasi-addi\-tive set function~$\Phi$.
  Lemma~\ref{lem: Pavlov} implies that
  $\Phi(U) = \Psi(U)$
  for all open
  $U \subset Y$.
\end{remark}


\section{Description of the Sobolev class in the~sense of Reshetnyak
for homeomorphisms}

This section is devoted to proving the following theorem.

\begin{theorem}\label{thm: comp oper Resh}
  Consider two metric measure spaces
  $(X, d_{X}, \mu)$
  and
  $(Y, d_{Y}, \nu)$,
  as well as a~homeomorphism
  $\varphi \colon X \to Y$.

  Then 
  $\varphi$
  induces the bounded composition operator
  $$
  \varphi^{*} \colon \Lip(Y) \to D^{1,q}(X), \quad 1 \le q < \infty,
  $$
  if and only if
  $\varphi \in D^{1,q}(X; Y)$
  in the sense of Reshetnyak.
  Furthermore,
  $\|\varphi^{*}\|^{q} = \int\limits_{X}|D_{q}\varphi|^{q}\, d \mu$.
  Moreover,
  the following set functions coincide:
  $$
  \| \varphi^{*} \vert U \|^{q}
  = \int\limits_{\varphi^{-1}(U)} |D_{q} \varphi|^{q}\, d \mu
  = \Psi(U)
  $$
  for every open set~$U$.
\end{theorem}

Before we proceed to proving Theorem~\ref{thm: comp oper Resh},
let us introduce some concepts.

Given a~family
$\mathcal{F}$
of functions in
$L^{p}(X, \mu)$,
a~function
$w\in L^{p}(X, \mu)$
with
$1 \le p < \infty$
is called a~\emph{majorant} for~$\mathcal{F}$
whenever
$|f| \le w$
a.e.\ for all
$f \in \mathcal{F}$.
The family~$\mathcal{F}$
is called~\emph{$o$-bounded} in
$L^{p}$
whenever at least one majorant of class
$L^{p}$
exists.
Denote by
$V_p(\mathcal{F})$
the collection of all class
$L^{p}$
majorants for~
$\mathcal{F}$.
Since
$L^{p}$
is a~$K$-space,
see~\cite[Ch.~VI (2.6)]{KVP},
the set
$V_p(\mathcal{F})$
contains the smallest element
$\sup \mathcal{F}$
in the sense of the natural order in
$L^{p}$.

We have the following~$o$-boundedness criterion
for sets in~$L^p$.

\begin{lemma}\label{lem: K lem}
  Given a~family
  $\mathscr{F} \subset L^{1}(X, \mu)$
  of nonnegative functions in
  $L^1(X, \mu)$,
  this family 
  is~$o$-bounded
  if and only if
  there exists a~constant
  $K < \infty$
  such that
  for every finite subset
  $\{f_1, f_2, \dots, f_k \} \subset \mathscr{F}$
  we have 
  \begin{equation}\label{eq: K lem eq}
    \int\limits_{X} \max \{ f_1, f_2, \dots, f_k \}\, d \mu \le K.
  \end{equation}
  Moreover,
  if~$\mathscr{F}$
  satisfies this condition
  then there exists a~sequence
  $(f_n)_{n \in \mathbb{N}}$
  of functions in~$\mathscr{F}$
  such that
  $$
  (\sup \mathscr{F})(x)
  = \lim\limits_{n \to \infty} \max\{f_1(x), f_2(x), \dots, f_n(x)\}
  $$
  for $\mu$-a.e.\ point
  $x \in X$.
\end{lemma}

Lemma~\ref{lem: K lem}
is a~particular case of the $o$-boundedness theorem
for sets in a~$K$-space with a~metric function,
see \cite[Ch. VI, Th.~1.31]{KVP}.
Repeating the ideas of \cite{KVP},
in the particular case
$\mathscr{F} \subset L^{1}(X, \mu)$
we give a~short proof of Lemma~\ref{lem: K lem}
in Appendix~\ref{sec: K lem},
see Lemma~\ref{lem: K lemm proof}.

Proceed to proving Theorem~\ref{thm: comp oper Resh}.

\subsection{Proof of necessity in Theorem~\ref{thm: comp oper Resh}}

Assume that~$\varphi$
induces the bounded composition operator.

Associated to the operator
$\varphi^{*}$,
there is the regular Borel measure~$\widetilde{\Psi}$,
see Theorem~\ref{thm: measure with lip}.
The definition of the measure~$\widetilde{\Psi}$
implies the following properties:

{\rm(1)}
$\widetilde{\Psi}(Y) =  \Psi(Y) = \|\varphi^{*}\|$;

{\rm(2)}
every function
$u \in \operatorname{1-Lip}(Y)$
and every Borel set
$U \subset Y$
satisfy 
\begin{equation}\label{eq:normative_ineq_q_infty}
  \int\limits_{\varphi^{-1}(U)}
  |\nabla_q(u\circ \varphi)|^{q}\, d \mu
  \le \Psi(U).
\end{equation}

Taking
$u, v \in \operatorname{1-Lip}(Y)$,
put
$$
U = \{ x\in X \mid |\nabla_q(u\circ \varphi)|(x)
\ge |\nabla_q(v \circ \varphi)|(x)\}
$$
and
$V = X \setminus U$.
Then~$U$ and~$V$
are disjoint Borel sets with
$U \cup V = X$.
Since~$\varphi$
is a~homeomorphism%
\footnote{%
  This is the only step in the proof
  relying on the property that~%
  $\varphi$
  is a~homeomorphism.
},
it follows that
$\varphi(U)$
and
$\varphi(V)$
are disjoint Borel sets with
$\varphi(U) \cup \varphi(V) = Y$.
Applying the inequality of claim~2 just above, 
we infer that
\begin{multline*}
  \int\limits_{X}
  \max\{|\nabla_{q} (u\circ \varphi)|^{q}, |\nabla_{q}(v\circ
  \varphi)|^{q}\}\, d \mu
  = \int\limits_{U} |\nabla_{q} (u\circ \varphi)|^{q}\, d \mu
  + \int\limits_{V} |\nabla_{q} (v\circ \varphi)|^{q}\, d \mu \\
  \le \Psi(\varphi(U)) + \Psi(\varphi(V))
  = \Psi(Y) = \|\varphi^{*}\|^{q}.
\end{multline*}
Similarly we can verify the inequality~\eqref{eq: K lem eq}
of Lemma~\ref{lem: K lem}
for a~finite tuple
$u_1$,
\dots,
$u_{k}$
of functions.
Lemma~\ref{lem: K lem} implies that
the family
$\{|\nabla_q(u\circ \varphi)| \mid u \in \operatorname{1-Lip}(Y)\}$
is~$o$-bounded in
$L_q(X)$.
Consequently,
$\varphi \in D^{1,q}(X; Y)$
in the sense of Reshetnyak.

From \eqref{eq:normative_ineq_q_infty}
we conclude that
$$
\int\limits_{\varphi^{-1}(U)} |D_{q} \varphi|^{q}\, d \mu \le \Psi(U)
= \|\varphi^{*} \vert U\|^{q}
$$
for every open set~$U$;
the last inequality holds by Lemma~\ref{lem: Pavlov}.

\subsection{Proof of sufficiency in Theorem~\ref{thm: comp oper Resh}}

Take some open set
$U \subset Y$
and
$u \in \operatorname{1-Lip}(Y)$.
Then
$u\circ \varphi \in D^{1,q}(X)$.
The chain rule yields
\begin{multline}
  \label{eq:suff_oper_Resh}
  \int\limits_{\varphi^{-1}(U)} |\nabla_q(u\circ \varphi)|^{q}\, d \mu
  \le
  \int\limits_{\varphi^{-1}(U)}
  |\nabla_q u|^{q}(\varphi(x)) |D_{q} \varphi|^{q} (x)
  \, d \mu(x)
  \\ \le
  \int\limits_{\varphi^{-1}(U)} |D_{q} \varphi|^{q}(x)\, d \mu(x).
\end{multline}
Inserting here
$U = Y$,
we see that
the composition operator
$\varphi^{*}$
is bounded.
In~addition,
\eqref{eq:suff_oper_Resh} implies that
$$
\Psi(U)
\le \int\limits_{\varphi^{-1}(U)} |D_{q} \varphi|^{q}\, d \mu
$$
for all open sets
$U \subset Y$.
\qed


\section{Proof of the main result}

\begin{theorem}\label{thm: main result}
  Consider two metric measure spaces
  $(X, d_{X}, \mu)$
  and
  $(Y, d_{Y}, \nu)$,
  where~$\mu$
  is a~Radon measure.
  Consider also a~homeomorphism
  $\varphi \colon X \to Y$
  and two numbers
  $1 \le q < p < \infty$
  with
  $\frac{1}{\sigma} = \frac{1}{q} - \frac{1}{p}$.

  The mapping~$\varphi$
  induces the bounded composition operator
  $$
  \varphi^{*} \colon D^{1,p}(Y) \cap \Lip(Y) \to D^{1,q}(X)
  $$
  if and only if
  the following properties hold:

  {\rm(1)}
  $\varphi \in D^{1,q}_{\operatorname{loc}}(X; Y)$;

  {\rm(2)}
  $\varphi$
  has finite distortion;

  {\rm(3)}
  the operator distortion function
  $K_{q,p}(x, \varphi)$
  lies in the space
  $L^{\sigma}(X)$.

  Furthermore,
  $\|K_{q,p}(\cdot, \varphi) \mid L^{\sigma}(X)\| = \|\varphi^{*}\|$.
  Moreover,
  the following set functions coincide:
  $$
  \Psi(U) =
  \int\limits_{\varphi^{-1}(U)} K_{q,p}(x, \varphi)^{\sigma}\, d \mu,
  $$
  for every open subset
  $U \subset Y$,
  where~$\Psi$
  is the measure of Theorem~{\rm\ref{thm: comp op measure}}.
\end{theorem}

\subsection{Proof of necessity in Theorem~\ref{thm: main result}}

Take a~homeomorphism
$\varphi \colon X \to Y$
which induces the bounded composition operator
$$
\varphi^{*} \colon D^{1,p}(Y) \cap \Lip(Y) \to D^{1,q}(X).
$$

If
$U \subset Y$
is an~open subset of~$Y$
with
$\nu(U) < \infty$,
for instance,
if~$U$
is a~metric ball of finite radius,
then every function
$u \in \Lip(Y) \cap D^{1,p}(Y)$
satisfies 
\begin{equation}\label{eq: bound for lip}
  \bigg(
    \int\limits_{\varphi^{-1}(U)}
    |\nabla_q (u\circ \varphi)|^{q}
    \, d \mu
  \bigg)^{\frac{1}{q}}
  \le
  \Psi(U)^{\frac{1}{\sigma}}
  \Lip(u; U) \nu(U)^{\frac{1}{p}}.
\end{equation}
From this we infer that~$\varphi$
induces the bounded composition operator
$$
\varphi^{*} \colon \Lip(U) \to D^{1,q}(\varphi^{-1}(U)).
$$
Consequently,
Theorem~\ref{thm: comp oper Resh} shows that
$\varphi \in D^{1,q}(\varphi^{-1}(U); U)$.
Since~$\varphi$
is continuous,
it follows that
$\varphi^{-1}(U)$
is open;
therefore,
each point
$x \in X$
has a~neighborhood~$W$
with
$\varphi \in D^{1,q}(W; Y)$.
From this we conclude that
$\varphi \in D^{1,q}_{\operatorname{loc}}(X; Y)$.
Moreover,
\eqref{eq: bound for lip} and Theorem~\ref{thm: comp oper Resh}
yield 
\begin{equation}\label{eq: bound for D}
  \bigg(
    \int\limits_{\varphi^{-1}(U)} |D_{q} \varphi|^{q}\, d \mu
  \bigg)^{\frac{1}{q}}
  \le
  \Psi(U)^{\frac{1}{\sigma}} \nu(U)^{\frac{1}{p}},
\end{equation}
where~$U$
is an~arbitrary open set of finite measure.

Verify that~$\varphi$
has finite distortion.

Denote by
$Z_{\varphi}$
the zero set of the Jacobian.
Then
$\nu(\varphi(Z_{\varphi})) = 0$.
Since~$\nu$
is a~regular measure,
for every
$\varepsilon > 0$
there exists an~open set~$U$
with 
$\varphi(Z_{\varphi}) \subset U$
and
$\nu(U) < \varepsilon$.
From~\eqref{eq: bound for D} we deduce that
$$
\bigg(
  \int\limits_{Z_{\varphi}} |D_{q} \varphi|^{q}\, d \mu
\bigg)^{\frac{1}{q}}
\le
\bigg(
  \int\limits_{\varphi^{-1}(U)} |D_{q} \varphi|^{q}\, d \mu
\bigg)^{\frac{1}{q}}
\le
\Psi(U)^{\frac{1}{\sigma}} \nu(U)^{\frac{1}{p}}
\le \|\varphi^{*}\| \varepsilon^{\frac{1}{p}}.
$$
Hence,
$|D_{q} \varphi| = 0$
is~$\mu$-a.e.\ in
$Z_{\varphi}$.

Let us verify that
$\widetilde{\Psi}$
is an~absolutely continuous measure with respect to~$\nu$.

Take some Borel set
$E \subset Y$
of zero~$\nu$-measure
and split it 
into three disjoint parts:
$E_1 = E \cap \Sigma_{\varphi^{-1}}$,
$E_2 = E \cap Z_{\varphi^{-1}}$,
and
$E_3 = E \setminus (E_1 \cup E_2)$.
It is clear that
$\widetilde{\Psi}(E_2 \cup E_3) = 0$
because in this case
$\mu(\varphi^{-1}(E_2 \cup E_3)) = 0$.
On the other hand,
$\varphi^{-1}(\Sigma_{\varphi^{-1}}) = Z_{\varphi}$
up to a~set of~$\mu$-measure zero,
while~$\varphi$
is a~mapping with finite distortion;
thus,
$\widetilde{\Psi}(E_2)=0$.

By the Radon--Nikodym theorem,
there exists a~function
$\Psi' \in L^{1}(Y, \nu)$
such that
$$
\int\limits_{U} \Psi' \, d \nu
= \widetilde{\Psi}(U)
$$
for every Borel set
$U \in \mathscr{B}(Y)$.

Verify that
$K_{q,p}(x, \varphi) \in L^{\sigma}(X, \mu)$.

From~\eqref{eq: bound for D} and absolute continuity of $\widetilde{\Psi}$
we find that
$$
\int\limits_{\varphi^{-1}(U)} |D_{q} \varphi|^{q}\, d \mu
\le
\Psi(U)^{\frac{q}{\sigma}} \nu(U)^{\frac{q}{p}} =
\bigg(
  \int\limits_{U}
  \Psi' \, d \nu
\bigg)^{\frac{q}{\sigma}}
\nu(U)^{\frac{q}{p}}
$$
for every open subset
$U \in Y$
of finite~$\nu$-measure.
Since~$\mu$ and~$\nu$
are regular measures
and~$\varphi$
is a~homeomorphism,
it follows that
the inequality
$$
\int\limits_{\varphi^{-1}(B)} |D_{q} \varphi|^{q}\, d \mu
\le
\bigg(
  \int\limits_{B}
  \Psi' \, d \nu
\bigg)^{\frac{q}{\sigma}}
\nu(B)^{\frac{q}{p}}
$$
holds for every Borel set
$B \in \mathscr{B}(Y)$.

On the other hand,
applying the change-of-variables formula in Lebesgue integral,
we find the chain of equalities
\begin{multline*}
  \int\limits_{\varphi^{-1}(U)} |D_{q} \varphi|^{q} \, d \mu =
  \int\limits_{\varphi^{-1}(U) \setminus Z_{\varphi}}
  \frac{|D_{q} \varphi|^{q}}{\mathcal{J}\varphi}
  \mathcal{J}(x,\varphi)\, d \mu\\
  =
  \int\limits_{U \setminus \varphi(Z_{\varphi} \cup \Sigma_{\varphi})}
  \frac{|D_{q} \varphi|^{q}\circ \varphi^{-1}}{\mathcal{J}\varphi
  \circ \varphi^{-1}} \, d \nu
  = \int\limits_{U} \chi_{Y \setminus \varphi(Z_{\varphi} \cup
  \Sigma_{\varphi})}
  \frac{|D_{q} \varphi|^{q}\circ \varphi^{-1}}{\mathcal{J}\varphi
  \circ \varphi^{-1}} \, d \nu.
\end{multline*}

Thus,
for every Borel set
$B \subset Y$
we have 
$$
\int\limits_{B}
\chi_{Y \setminus \varphi(Z_{\varphi} \cup \Sigma_{\varphi})}
\frac{|D_{q} \varphi|^{q}\circ \varphi^{-1}}
{\mathcal{J}\varphi \circ \varphi^{-1}}
\, d \nu
\le
\bigg(
  \int\limits_{B}
  \Psi' \, d \nu
\bigg)^{\frac{q}{\sigma}}
\nu(B)^{\frac{q}{p}}.
$$
Applying the ``reverse'' H\"older inequality,
see Appendix~\ref{sec: reversed Holder},
and Proposition~\ref{prp: reversed Holder},
we obtain the inequality
$$
\frac{|D_{q} \varphi|^{q}(\varphi^{-1}(y))}{\mathcal{J}
\varphi(\varphi^{-1}(y))}
\le \big(\Psi'(y)\big)^{\frac{q}{\sigma}}
\quad
\text{for $\nu$-a.e.\ }
\quad
y \in Y \setminus \varphi(Z_{\varphi} \cup \Sigma_{\varphi}).
$$

Outside the set
$\varphi(Z_{\varphi}) \cup \varphi(\Sigma_{\varphi})$
the mapping
$\varphi^{-1}$
has Luzin's
$\mathcal{N}$-property.
Thus,
we have 
$$
\frac{|D_{q} \varphi|^{q}(x)}{\mathcal{J} \varphi(x)^{\frac{1}{p}}}
\le
\big(
  \Psi'(\varphi(x)) \mathcal{J}\varphi(x)
\big)^{\frac{1}{\sigma}}
\quad
\text{for $\mu$-a.e.\ }
\quad
x \in X \setminus (Z_{\varphi} \cup \Sigma_{\varphi}).
$$
Since the set
$\Sigma_{\varphi}$
is of~$\mu$-measure zero
and the operator distortion function
$K_{q,p}(x, \varphi)$
vanishes for
$x \in Z_{\varphi}$,
we arrive at 
\begin{equation}\label{eq: final estim}
  K_{q,p}(x, \varphi)
  \le
  \big(
    \Psi'(\varphi(x)) \mathcal{J}(x)
  \big)^{\frac{1}{\sigma}}
  \quad
  \text{for $\mu$-a.e.\ }
  \quad
  x \in X.
\end{equation}

This inequality 
yields the estimate
$$
\int\limits_{\varphi^{-1}(U)}
K_{q,p}(x,\varphi)^{\sigma}
\, d \mu(x)
\le
\int\limits_{U} \Psi'(y)\, dy
= \Psi(U)
$$
for every open set
$U \subset Y$.
\qed

\subsection{Proof of sufficiency in Theorem~\ref{thm: main result}}

Take
$u \in D^{1,p}(Y) \cap \Lip(Y)$.
The chain rule yields
$u\circ \varphi\in D^{1,q}_{\operatorname{loc}}(X)$.
Denote by
$g_{u\circ \varphi}$
an~upper gradient for the function
$u\circ \varphi$
and by
$g_{u}$
an~upper gradient for~$u$.
Applying the chain rule,
we then infer that
$$
g_{u\circ \varphi}(x) \le g_{u}(\varphi(x)) |D_{q} \varphi|(x)
$$
for $\mu$-a.e.\
$x \in X$.
Furthermore,
since the distortion of~$\varphi$ is finite,
H\"older's inequality
($\frac{q}{\sigma} + \frac{q}{p} = 1$)
and the change-of-variables formula in Lebesgue integral
imply that
\begin{multline*}
  \int\limits_{X} g_{u\circ \varphi}^{q}(x) \, d \mu (x)
  \le
  \int\limits_{X} g_{u}^{q}(\varphi(x)) |D_{q} \varphi|^{q}(x)\, d \mu(x)
  \\
  =
  \int\limits_{X \setminus Z_{\varphi}}
  g_{u}^{q}(\varphi(x))\mathcal{J} \varphi(x)^{\frac{q}{p}} \cdot
  \frac{|D_{q} \varphi|^{q}(x)}{\mathcal{J} \varphi(x)^{\frac{q}{p}}}
  \, d \mu(x) \\
  \le
  \bigg(
    \int\limits_{X \setminus Z_{\varphi}}
    g_{u}^{p}(\varphi(x)) \mathcal{J} \varphi(x)
    \, d \mu(x)
  \bigg)^{\frac{q}{p}}
  \bigg(
    \int\limits_{X \setminus Z_{\varphi}}
    \bigg(\frac{|D_{q} \varphi|(x)}{\mathcal{J}
    \varphi(x)^{\frac{1}{p}}} \bigg)^{\sigma}
    \, d \mu(x)
  \bigg)^{\frac{q}{\sigma}} \\
  \le
  \bigg(
    \int\limits_{Y}
    g_{u}^{p}(y)\, d \nu(y)
  \bigg)^{\frac{q}{p}}
  \bigg(
    \int\limits_{X}
    K_{q,p}(x, \varphi)^{\sigma}
    \, d \mu(x)
  \bigg)^{\frac{q}{\sigma}}.
\end{multline*}
Hence,
$u\circ \varphi \in D^{1,q}(X)$.

There exists a~sequence
$(g_{n})_{n\in \mathbb{N}}$
of upper gradients for~$u$
such that
$g_n \to |\nabla_p u|$
in
$L^{p}(Y, \nu)$
as
$n \to \infty$,
see \cite[Theorem 6.3.20]{HKNT15} for instance.
Then 
\begin{multline*}
  \bigg(
    \int\limits_{X}
    |\nabla_{q} (u\circ \varphi)|^{q}(x)\, d \mu (x)
  \bigg)^{\frac{1}{q}}
  \le
  \bigg(
    \int\limits_{X} g_{u\circ \varphi}^{q}(x) \, d \mu (x)
  \bigg)^{\frac{1}{q}}\\
  \le
  \bigg(
    \int\limits_{X}
    K_{q,p}(x, \varphi)^{\sigma}
    \, d \mu(x)
  \bigg)^{\frac{q}{\sigma}}
  \bigg(
    \int\limits_{Y}
    g_{n}^{p}(y)\, d \nu(y)
  \bigg)^{\frac{q}{p}}
  \\
  \to
  \bigg(
    \int\limits_{X}
    K_{q,p}(x, \varphi)^{\sigma}
    \, d \mu(x)
  \bigg)^{\frac{q}{\sigma}}
  \bigg(
    \int\limits_{Y}
    |\nabla_p u|^{p}(y)\, d \nu(y)
  \bigg)^{\frac{q}{p}}
  \quad \text{as} \quad n \to \infty.
\end{multline*}

In~particular,
with the replacement of~$X$
by~$\varphi^{-1}(U)$
and~$Y$
by~$U$
in the previous calculations,
this implies that
every open subset
$U \subset Y$
satisfies 
$$
\Psi(U) \le \int\limits_{\varphi^{-1}(U)} K_{q,p}(x,
\varphi)^{\sigma}\, d \mu(x).
$$

Thus,
the proof of Theorem~\ref{thm: main result} is complete.
\qed

\begin{remark}
  Proposition~\ref{prp: funct eq}
  shows that
  the following set functions coincide:
  $$
  \|\varphi^{*} \vert U \|^{\sigma}
  = \Phi(U) = \Psi(U)
  = \int\limits_{\varphi^{-1}(U)} K_{q,p}(x, \varphi)^{\sigma}\, d \mu(x)
  $$
  for all open sets
  $U \subset Y$
  of finite measure,
  where~$\Phi$
  is the~quasi-additive set function,
  see Section~\ref{par: set func}.

  Therefore,
  the function~$\Phi$
  is the trace of the measure~$\Psi$
  on the open sets of finite measure.
\end{remark}

By analogy with the authors' previous article~\cite{VS25},
the proof of Theorem~\ref{thm: main result}
uses the following properties of the space
$D^{1,p}$:

{\rm(1)}
the space
$D^{1,p}$
is closed under the pointwise maximum and minimum of functions,
and furthermore,
$$
|\nabla_p (u_1 \vee u_2)|
= |\nabla_p u_1| \cdot \chi_{\{u_1 \ge u_2\}}
+ |\nabla_p u_2| \cdot \chi_{\{u_2 > u_1\}}, \quad u_1, u_2 \in D^{1,p},
$$
almost everywhere.
The similar property of
$u_1 \wedge u_2$
is also required.

{\rm(2)}
Lipschitz functions lie in the space
$D^{1,p}_{\operatorname{loc}}$,
and furthermore,
$|\nabla_p u|(y) \le \operatorname{Lip} u (y)$
almost everywhere for all
$u \in \operatorname{Lip} \cap D^{1,p}$,
where
$$
\operatorname{Lip}u (y) =
\varlimsup\limits_{r \to 0}
\frac{\sup\{|u(y')-u(y)| \mid d_{Y}(y',y)\le r\}}{r}.
$$

{\rm(3)}
Given a~measurable mapping
$\varphi \colon X \to Y$
of class
$D^{1,p}(X;Y)$
in the sense of Reshetnyak
and a~Lipschitz function
$u \in \operatorname{Lip}(Y)$,
the composition
$u \circ \varphi$
lies in 
$D^{1,p}(X)$
and the ``norm of the gradient'' satisfies
$$
|\nabla_p(u\circ \varphi)|(x) \le |\nabla_p u|(\varphi(x))|D_p \varphi|(x)
\quad \text{almost everywhere.}
$$

Therefore,
if the Sobolev space
$S^{p}$,
where
$1 \le p < \infty$,
whose functions are defined on a~metric measure space,
has properties (1)--(3)
then Theorem~\ref{thm: main result} holds.

\begin{theorem}
  Consider two metric measure spaces
  $(X, d_{X}, \mu)$
  and
  $(Y, d_{Y}, \nu)$,
  where~$\mu$
  is a~Radon measure.
  A~homeomorphism
  $\varphi \colon X \to Y$
  induces the bounded composition operator
  $$
  \varphi^{* } \colon S^{p}(Y) \cap \operatorname{Lip}(Y) \to S^{q}(X),
  \quad 1 \le q \le p < \infty,
  $$
  if and only if
  the following properties hold:

  {\rm(i)}
  $\varphi \in S^{q}_{\operatorname{loc}}(X; Y)$
  in the sense of Reshetnyak;

  {\rm(ii)}
  $\varphi$
  has finite distortion;

  {\rm(iii)}
  the operator distortion function
  $K_{q,p}(\cdot, \varphi)$
  lies in the space
  $L^{\sigma}(X, \mu)$,
  where
  $\frac{1}{\sigma} = \frac{1}{q} - \frac{1}{p}$.

  The norm of the composition operator
  $\|\varphi^{*}\|$
  coincides with
  the norm of $K_{q,p}$ in $L^{\sigma}(X, \mu)$.
  Moreover,
  for
  $q < p$
  the following set functions coincide:
  $$
  \Psi(U) =
  \int\limits_{\varphi^{-1}(U)} K_{q,p}(x, \varphi)^{\sigma}\, d \mu,
  $$
  for every open
  $U \subset Y$.
\end{theorem}


\appendix

\section{Proof of the ``reverse'' H\"older's inequality}
\label{sec: reversed Holder}

\begin{propos}
  \label{prp: reversed Holder}
  Consider three positive measurable functions~$f$,
  $g$ and~$h$
  on some measure space
  $(Y, \Sigma, \nu)$
  and take two numbers
  $p$,
  $q \ge 1$
  with
  $\frac{1}{p} + \frac{1}{q} = 1$.

  The inequality
  $$
  \int\limits_{U} f(y) \, d \nu(y)
  \le
  \bigg(\int\limits_{U} g(y) \, d \nu(y) \bigg)^{\frac{1}{p}}
  \bigg(\int\limits_{U} h(y) \, d \nu(y) \bigg)^{\frac{1}{q}}
  $$
  holds for every
  $U \in \Sigma$
  if and only if
  $f(y) \le g(y)h(y)$
  holds $\nu$-a.e.
\end{propos}

\begin{proof}
  Sufficiency is quite obvious.
  Let us verify necessity.

  Assume on the contrary that
  $\{ y \in Y \mid f(y) > g(y)h(y)\}$
  is of positive measure.
  We may assume that
  the functions
  $g(y)$
  and
  $h(y)$
  are not vanishing at a.e.\ point
  $y\in Y$.
  Then the set
  $
  \{y\in Y \mid f(y) \ge r+\varepsilon, \ g(y)h(y) \le r \}
  $
  is of positive measure for some
  $r > 0$
  and
  $\varepsilon > 0$.
  Furthermore,
  the set
  $$
  \{y\in Y \mid f(y) \ge r + \varepsilon,\ g(y) \le a,\ g(y)h(y) \le r \}
  $$
  is of positive measure for some
  $a > 0$.
  Choose now
  $n_0 \ge 2$
  such that
  the set
  $$
  \{y\in Y \mid f(y) \ge r + \varepsilon,\ \frac{a}{n_0} \le g(y) \le
  a,\ g(y)h(y) \le r \}
  $$
  is of positive measure.

  Consider the sets
  $$
  E^{1}_{j} = \{y\in Y \mid f(y) \ge r + \varepsilon,\
  a\frac{j}{n_0} \le g(y) \le a\frac{j+1}{n_0},\ g(y)h(y) \le r \},
  $$
  where
  $j=1,2,\dots, n_0-1$.
  Then there is some index
  $j_1 \ge 1$
  such that
  the set
  $E^{1}_{j_1}$
  is of positive measure.
  Denote
  $E^{1}_{j_1}$
  by
  $E^{1}$.

  Suppose that
  $k \ge 2$
  and that
  the set
  $E^{k-1} = E^{k-1}_{j_{k-1}}$
  is chosen.
  Consider the set
  $$
  E^{k}_{j} = \{y \in E^{k-1} \mid f(y) \ge r+ \varepsilon, \
  \frac{j}{n_0^{k}} \le g(y) \le \frac{j+1}{n_0^{k}}, \ g(y)h(y) \le r\},
  $$
  where
  $j=1,2,\dots, n_0^{k}-1$.
  Since
  $\nu(E^{k-1}) > 0$,
  there is
  $j_k$
  such that
  $\nu(E^{k}_{j_{k}}) > 0$;
  furthermore,
  $j_{k} \ge j_{k-1} n_0$
  because
  $E^{k}_{j_k} \subset E^{k-1} = E^{k-1}_{j_{k-1}}$.
  Put
  $E^{k} = E^{k}_{j_k}$.

  Therefore,
  we obtain two sequences of positive integers
  $j_k$
  and sets
  $E^{k}$
  such that:

  {\rm(1)}
  $(j_{k})_{k\in \mathbb{N}}$
  is an~increasing sequence
  and 
  $n_0^{k} \le j_{k}$;

  {\rm(2)}
  for
  $y\in E^{k}$
  we have 
  $$
  f(y) \ge r + \varepsilon,
  \quad
  a \frac{j_{k}}{n_0^{k}} \le g(y) \le a \frac{j_{k}+1}{n_0^{k}},
  \quad
  g(y)h(y) \le r,
  $$
  which implies that
  $h(y) \le \frac{n_0^{k} r}{a j_k}$
  for
  $y \in E^{k}$;

  {\rm(3)}
  $\nu(E^{k}) > 0$
  for all
  $k \in \mathbb{N}$.

  Therefore,
  \begin{align*}
    \nu(E^{k}) (r + \varepsilon)
    & \le
    \int_{E^{k}} f\, d \nu
    \le
    \bigg(
      \int\limits_{E_k}
      g^{p} \, d \nu
    \bigg)^{\frac{1}{p}}
    \bigg(
      \int\limits_{E_k}
      h^{q} \, d \nu
    \bigg)^{\frac{1}{q}}
    \\
    & \le a \frac{j_{k}+1}{n_{0}^{k}} \nu(E^{k})^{\frac{1}{p}}
    \cdot \frac{n_0^{k} r}{a j_k} \nu(E^{k})^{\frac{1}{q}}
    = \nu(E^k)\frac{j_{k}+1}{j_k} r.
  \end{align*}
  From this we infer that
  $r+ \varepsilon \le r \frac{j_{k}+1}{j_k}$,
  but
  $j_k \to \infty$
  as
  $k \to \infty$,
  and so
  $r + \varepsilon \le r$,
  which is a~contradiction.
\end{proof}

\section{Proof of the lemma about $o$-boundedness}
\label{sec: K lem}

\begin{lemma}
  \label{lem: K lemm proof}
  A~family
  $\mathscr{F} \subset L^{1}(X, \mu)$
  of nonnegative functions 
  is~$o$-bounded
  if and only if
  for some constant
  $K < \infty$
  every finite subset
  $\{f_1, f_2, \dots, f_k \} \subset \mathscr{F}$
  satisfies 
  $$
  \int\limits_{X} \max \{ f_1, f_2, \dots, f_k \}\, d \mu \le K.
  $$
  Moreover,
  if~$\mathscr{F}$
  meets this condition
  then there exists a~sequence
  $(f_n)_{n \in \mathbb{N}}$
  of functions in~$\mathscr{F}$
  such that
  $$
  (\sup \mathscr{F})(x)
  = \lim\limits_{n \to \infty} \max\{f_1(x), f_2(x), \dots, f_n(x)\}
  $$
  for $\mu$-a.e.\ point
  $x \in X$.
\end{lemma}

\begin{proof}
  Necessity is obvious.
  Let us verify sufficiency.

  We may assume that
  the system~$\mathscr{F}$
  is closed under maxima.
  By the hypotheses of the lemma,
  every function
  $f \in \mathscr{F}$
  satisfies 
  $$
  \int\limits_{X} f \, d \mu \le K.
  $$
  Hence,
  $
  \sup\big\{\int\limits_{X} f\, d \mu \mathrel{\big\vert} f \in
  \mathscr{F}\big\} = M
  $
  is finite
  and there exists a~sequence of functions
  $(\widetilde{f}_n)_{n \in \mathbb{N}} \subset \mathscr{F}$
  with
  $$
  M = \lim\limits_{n \to \infty} \int\limits_{X} \widetilde{f}_{n} \, d \mu.
  $$
  Consider the monotone sequence
  $$
  f_{n}
  = \max \{ \widetilde{f}_1, \widetilde{f}_2, \dots, \widetilde{f}_{k}\}
  \in \mathscr{F}.
  $$
  Then
  $f_n$
  increases monotonely
  and the integrals
  $\int\limits_{X} f_n\, d \mu$
  are bounded by the constant~$K$.
  Consequently,
  by Beppo Levi's theorem the limit function
  $w = \lim\limits_{n \to \infty} f_n$
  lies in 
  $L^1(X)$
  and 
  \begin{equation}\label{eq:this_is_sup}
    M = \int\limits_{X} w\, d \mu.
  \end{equation}
  Since the functions
  $f \in \mathscr{F}$
  are positive,
  we conclude that~$w$
  is a~majorant for the family~$\mathscr{F}$.
  By~\eqref{eq:this_is_sup}
  the function~$w$
  coincides with
  $\sup \mathscr{F}$
  for $\mu$-a.e.\ point of~$X$.
\end{proof}

\section*{Acknowledgments} 
The work of D.~A.~Sboev
was supported by the Mathematical Center in Aka\-dem\-gorodok
under Agreement No.~075--15--2025--349
with the Ministry of Science and Higher Education of the Russian Federation.

The work of S.~K.~Vodopyanov was carried out
in the framework of the State Task to the Sobolev Institute of Mathematics
(Project FWNF--2022--0006).

\printbibliography

\end{document}